\documentclass[10pt,twosided]{article}
\usepackage[tbtags]{amsmath}
\usepackage{amssymb}
\allowdisplaybreaks[4]
\usepackage{color}

\usepackage{amssymb,color}

\definecolor{c20}{rgb}{0.,0.7,0.}
\definecolor{c30}{rgb}{0.,0.,1.}
\definecolor{c40}{rgb}{1,0.1,0.7}
\definecolor{c50}{rgb}{1,0,0}
\definecolor{c60}{rgb}{0,0.9,0.1}

\def\aH#1{\textcolor{c20}{#1}}
\def\aH#1{#1}

\def\rH#1{\textcolor{c20}{#1}}
\def\rH#1{#1}

\def\Hr#1{\textcolor{c20}{#1}}
\def\Hr#1{#1}

\def\xH#1{\textcolor{c20}{#1}}
\def\xH#1{#1}
\def\xx#1{\textcolor{c60}{#1}}
\def\xx#1{#1}

\def\cJ#1{\textcolor{c50}{#1}}
\def\cJ#1{#1}

\def\cJI#1{\textcolor{c50}{#1}}
\def\cJI#1{#1}
\def\cL#1{\textcolor{c50}{#1}}
\def\cL#1{#1}
\def\ccL#1{\textcolor{c50}{#1}}
\def\ccL#1{#1}

\def\rL#1{\textcolor{c50}{#1}}
\def\rL#1{#1}

\def\rP#1{\textcolor{c50}{#1}}
\def\rP#1{#1}

\usepackage{amssymb,color}

\newcommand{\abs}[1]{\lvert #1 \rvert}
\newcommand{\ABs}[1]{ \biggl \lvert #1 \biggr \rvert}

\newcommand{\E}[1]{\mathbb{E}\left(#1\right)}

\newcommand{\pk}[1]{\mathbb{P} \left\{ #1 \right\} }

\newcommand{\R}{\mathbb{R}}

\newcommand{\inr}{\in \R}

\newcommand{\limit}[1]{\lim_{#1 \to   \infty}}
\newcommand{\todis}{\stackrel{d}{\to}}

\def\Ha{ \mathcal{H}_{\alpha}}

\def\CC{\mathbb{C}}

\newcommand{\BQN}{\begin{eqnarray}}
\newcommand{\EQN}{\end{eqnarray}}
\newcommand{\BQNY}{\begin{eqnarray*}}
\newcommand{\EQNY}{\end{eqnarray*}}

\newcommand{\BS}{\begin{sat}}
\newcommand{\ES}{\end{sat}}
\newcommand{\BT}{\begin{theo}}
\newcommand{\ET}{\end{theo}}
\newcommand{\BK}{\begin{korr}}
\newcommand{\EK}{\end{korr}}

\newcommand{\BD}{\begin{de}}
\newcommand{\ED}{\end{de}}

\newcommand{\BIT}{\begin{itemize}}
\newcommand{\EIT}{\end{itemize}}
\newcommand{\BDI}{\begin{description}}
\newcommand{\EDI}{\end{description}}

\newcommand{\BRM}{\begin{remarks}}
\newcommand{\ERM}{\end{remarks}}

\newcommand{\BEL}{\begin{lem}}
\newcommand{\EEL}{\end{lem}}

\newtheorem{theo}{Theorem}[section]
\newtheorem{sat}[theo]{Proposition}
\newtheorem{de}[theo]{Definition}
\newtheorem{lem}[theo]{Lemma}

\newtheorem{korr}[theo]{Corollary}
\newtheorem{remark}[theo]{Remark}
\newtheorem{remarks}[theo]{Remarks}

\newcommand{\nelem}[1]{{Lemma \ref{#1}}}

\newcommand{\netheo}[1]{{Theorem \ref{#1}}}

\newcommand{\prooftheo}[1]{ \textsc{Proof of Theorem} \ref{#1} }

\newcommand{\COM}[1]{}

\newcommand{\QED}{\hfill $\Box$}

\def\Ga{\gamma}

\def\rw{\rightarrow}

\def\IF{\infty}
\def\piter{\mathcal{P}}

\def\tildet{\tilde{t}_0}

\date{}

\def\mlH{\mathcal{H}}

\def\oo{(1+o(1))}
\def\asu{\ \ \text{as}\ u\rw\IF}
\def\tauu{\tau_1(u)}
\def\tauuu{\tau_2(u)}
\def\tAu{\tau_1^*(u)}
\def\tAuu{\tau_2^*(u)}
\newcommand{\equaldis}{\stackrel{d}{=}}

\def\deluu{\delta_2(u)}
\def\delu{\delta_1(u)}
\def\Dux{\widetilde{\Delta_x^{1}}(u)}
\def\Duux{\widetilde{\Delta_x^{2}}(u)}
\def\tilu{\tilde{u}}
\def\tilt{\tilde{t}_0}
\def\CC{\mathbb{C}}
\def\LT{\left}
\def\RT{\right}
\def\FF{\mathcal{F}}
\def\FFab{\FF_{\alpha,\beta}}
\def\det{\triangle}
\def\wdet{\widetilde{\triangle}}
\def\H{\mathcal{H}}

\newcommand{\toprob}{ \stackrel{p}{\to}}
\begin{document}

\title{\bf Approximation \rH{of Passage Times} of $\gamma$-reflected Processes \rH{with fBm Input}}

\bigskip
\author{ Enkelejd Hashorva and Lanpeng Ji\thanks{University of Lausanne, UNIL-Dorigny 1015 Lausanne, Switzerland, enkelejd.hashorva@unil.ch, lanpeng.ji@unil.ch}
}

 \maketitle
\vskip -0.61 cm


\bigskip
{\bf Abstract:} Define \cJI{a} $\gamma$-reflected process $W_\Ga(t)=Y_H(t)-\Ga\inf_{s\in[0,t]}Y_H(s)$, $t\ge0$ with input process $\{Y_H(t), t\ge 0\}$ \aH{which is} a fractional Brownian motion with Hurst index $H\in (0,1)$ and a negative linear trend. In risk theory
$R_\gamma(u)=u-W_\Ga(t),  t\ge0$ \rH{is referred to as} the risk process with tax payments of a {loss-carry-forward} type.
\aH{For various risk processes numerous results are known for the approximation of the first \aH{and last} passage times to 0 (ruin times) when the initial reserve  $u$ goes to infinity.  In this paper we show that for the $\gamma$-reflected process the conditional (standardized) first and last passage times
are jointly asymptotically Gaussian and completely dependent. An important contribution of this paper is that it links ruin problems with
extremes of non-homogeneous Gaussian random fields defined by $Y_H$ which are also investigated in this contribution.}

{\bf Key Words:}   \rH{Gaussian} approximation; passage times; \rH{$\gamma$-reflected process};   workload process; \Hr{risk process with tax}; fractional Brownian motion; \rH{Piterbarg constant; Pickands constant}.\\

{\bf AMS Classification:} Primary 60G15; secondary 60G70

\section{Introduction and Main Result}
Let $\{X_H(t), t\ge0\}$ be  a standard fractional Brownian motion (fBm) with Hurst index  $H\in(0,1)$ meaning that $X_H$ is
a  centered Gaussian process with covariance  \xH{function}
\BQNY
Cov(X_H(t),X_H(s))=\frac{1}{2}(t^{2H}+s^{2H}-\mid t-s\mid^{2H}),\quad t,s\ge0.
\EQNY
We \rH{shall define} the $\gamma$-reflected process with input process $Y_H(t)= X_H(t)- ct$ by
\BQN\label{Wgam}
W_\Ga(t)=Y_H(t)-\Ga\inf_{s\in[0,t]}Y_H(s), \  \ t\ge0, \label{PCW}
\EQN
where $\gamma \in [0,1]$ and $c>0$ are two fixed constants.\\
\rH{Motivations for studying $W_\Ga$ come from both risk and queuing theory. For instance, in queuing theory} $W_1$ is the so-called  workload process (or queue length process), see e.g., Harrison (1985), Asmussen (1987), Zeevi and Glynn (2000), Whitt (2002) and Awad and Glynn (2009) among many others.
In advanced risk  theory the process  $R_\gamma(t)=u-W_\Ga(t),  t\ge0,u\ge 0$ is referred to as
the risk process with tax payments of a loss-carry-forward type, see e.g., Asmussen and Albrecher (2010).\\
\rH{Recently Hashorva et al.\ (2013) studed the \ccL{asymptotics of the} probability $\pk{\sup_{t\in [0,T]} W_\gamma(t)> u}$
as $u\to \IF$ for \ccL{both $T<\IF$ and}  $T=\IF$.  Continuing the investigation of the aforementioned paper
in this contribution we shall investigate the approximation of first and last passage times of \ccL{$W_\Ga$}}. Specifically,  define the first and last  passage times of $W_\Ga$ to a constant threshold $u>0$ by
\BQN\label{deftau}
\tau_1(u)=\inf\{t\ge0, W_\Ga(t)>u\}\ \ \ \text{and}\ \ \ \tau_2(u)=\sup\{t\ge0, W_\Ga(t)>u\},
\EQN
respectively (here we use that $\inf\{\emptyset\}=\IF$). Further,
define $\tAu,\  \tAuu, u>0$ in the same probability space
such that
\BQN\label{eq:taustar}
(\tAu, \tAuu) \equaldis  (\tauu,\tauuu)   \Bigl\lvert (\tauu < \IF),
\EQN
where $\equaldis$ stands for equality of distribution functions.\\
The first and last passage times of Gaussian processes conditioned \ccL{on} that $\tauu < \IF$ are analysed in H\"usler and Piterbarg (2008)
and H\"usler and Zhang (2008) \ccL{when $\Ga=0$}.  Therein, the Guassian approximations of both $\tAu$ and $\tAuu$ are \rH{derived} as $u\to \IF$.
The Gaussian approximation is \ccL{not only} of theoretical interest \ccL{but also} important for statistical estimation.
First passage times (sometimes called ruin times) are also studied extensively in the framework of insurance risk processes,
see the recent articles Griffin and Maller (2012), Griffin (2013), Griffin et al. (2013), D\c{e}bicki et al.\ (2013)
and the monographs Embrechts et al.\ (1997), Asmussen and Albrecher (2010) \rH{for approximations} of
ruin times of various risk processes. In this framework, $\tAu$ can be \rH{interpreted} as the conditional ruin time of the fBm risk process with tax payments of a {loss-carry-forward} type.\\
With motivation from the aforementioned \rH{contributions}, this paper is concerned with the Gaussian approximation of the random vector  $(\tAu,
\tAuu)$, as $u \to \IF$. \rH{For the derivation of the tail asymptotics of $\sup_{t\in [0,T]} W_\gamma(t)$
Hashorva et al.\ (2013) showed that the investigation of the supremum of certain \rL{non-stationary} Gaussian random fields is crucial.
One key merit of our problem of approximating the joint distribution function of $(\tAu,
\tAuu)$ is that it \ccL{leads}, as in the case of  the analysis of the tail asymptotics of  $\sup_{t\in [0,T]} W_\gamma(t)$,
to \Hr{an} interesting unsolved problem of asymptotic theory of Gaussian random fields. Although the latter investigation
 was not initially in the scope of this paper, the result derived in \ccL{\netheo{Thmx} is} important for various theoretical questions.}
\rH{Next, set}
%
%
$$
A(u)=\frac{H^{H+1/2}}{(1-H)^{H+1/2}c^{H+1}}u^H,\ \ \ \text{and}\ \ \ \tilde{t}_0=\frac{H}{c(1-H)}
$$
and denote by $\todis$ and $\toprob$ the convergence in distribution and in probability, respectively.
Further, let \ccL{$\mathcal{N}$} be a  $N(0,1)$ random variable. Our \rH{principal result is} the following theorem:
\BT\label{Thmfl}
Let the $\gamma$-reflected process $\{W_\Ga(t),t\ge0\}$ be given as in \eqref{Wgam} with $\Ga\in(0,1)$, and let $\tAu, \tAuu$ be defined as in \eqref{eq:taustar}.
Then,  as $u\rw\IF$
\BQN\label{eq:asym}
\Biggl( \frac{\tAu-\tilde{t}_0u}{A(u)} , \frac{\tAuu-\tilde{t}_0u}{A(u)}  \Biggr) &\todis &(\mathcal{N}, \mathcal{N}).
\EQN

\ET

{\bf Remarks}: a) 
The joint convergence in \eqref{eq:asym} implies $({\tAuu-\tAu})/{A(u)} \toprob 0$ \rH{as $u\to \IF$}.\\
b) For any $u\ge0$
$\pk{\tauu<\IF}=1$  when $\Ga=1$ (cf. Duncan and Jin (2008)),  which is the reason \ccL{of}
considering only the case that $\Ga\in(0,1)$. Under the latter assumption \aH{on $\gamma$} we have further that 
$\pk{\tauuu<\IF\bigl\lvert \tauu<\IF}\Hr{=1}$, \Hr{which
follows from the} fact that  $\lim_{t\rw\IF}W_\Ga(t)=-\IF$ almost surely \Hr{since in view of Remark 5 in Kozachenko et al.\ (2011)}
$$
\lim_{t\rw\IF}\frac{\sup_{s\in[0,t]}\abs{X_H(s)}}{t}=0, \quad \forall H\in(0,1).
$$

c) It is surprising that the Gaussian approximation of the conditional first and last passage times does not involve the reflection constant $\gamma$. \\

Organisation of the rest of the paper: In the next section we present a key result on \aH{the} supremum of \aH{some} Gaussian random fields \aH{defined} by  $Y_H$ and then display the proof of \netheo{Thmfl}. Section 3 is dedicated to the proof of Theorem 2.1. A variant of Piterbarg Lemma suitable for 
 Gaussian random fields is presented in Appendix.

\section{Further Results and Proof of \netheo{Thmfl}}
Following the idea of H\"{u}sler and Piterbarg (1999, 2008), and as discussed in Hashorva et al. (2013) it is convenient to
introduce the following family of Gaussian random fields:
$$Y_u(s,t):=\frac{X_H(ut)-\gamma X_H(us)}{(1+ct-c\Ga s)u^H},\ \  s,t\ge0.$$
The variance function of $\{Y_u(s,t), s,t\ge0\}$ is given by
\BQN\label{eq:varY}
V_Y^2(s,t)=\frac{(1-\gamma)t^{2H}+(\Ga^2-\Ga)s^{2H}+\Ga(t-s)^{2H}}{(1+ct-c\Ga s)^2},\ \   s,t\ge0.
\EQN
Moreover,  on the set $\{(s,t): 0\le s\le t<\IF\}$ it attains its maximum at  the unique point $(0,\tilt)$
with $\tilde{t}_0=\frac{H}{c(1-H)}$ and further
$$
V_Y(0, \tilde{t}_0)=\frac{H^H (1-H)^{1-H}}{c^H }.
$$

By changing time $t=t'u, s=s'u$ and noting that the distribution of $Y_u$ does not depend on $u$, we obtain
\BQNY 
\pk{\tauu<\IF}&=&\pk{\exists t\in[0,\IF) \ \text{such that}\ W_\Ga(t)>u}\nonumber\\
&=&\pk{\exists  t' \in[0,\IF)\ \text{such that}\ Y_u(s',t')>u^{1-H} \ \text{for some} \  s'\in[0, t']}\nonumber\\
&=&\pk{\exists  t \in[0,\IF)\ \text{such that}\ Y(s,t)>u^{1-H} \ \text{for some} \  s\in[0, t]},
\EQNY
where
\BQN\label{Yst}
Y(s,t):=\frac{X_H(t)-\gamma X_H(s)}{1+\rH{c(t-\Ga s)}},\ \  s,t\ge0.
\EQN
\COM{
Define
\BQNY
T_1(u)=\inf\{t\ge0:\ Y(s,t)>u^{1-H}\ \text{for some} \  s\in[0, t]\}\ 
\EQNY
and
\BQNY
 T_2(u)=\sup\{t\ge0: Y(s,t)>u^{1-H}\ \text{for some} \  s\in[0, t]\}.
\EQNY
Clearly $\tau_i(u)=u T_i(u), i=1,2$. Further
for any $x\inr$ 
\BQN\label{eq:tau1}
\pk{\frac{\tau_1(u)-\tilde{t}_0u}{A(u)}\le x\Big| \tauu<\IF}
&=&\frac{\pk{\sup_{0\le s\le t\le\tilde{t}_0  +xA(u)u^{-1}}Y(s,t)>u^{1-H}}}{\pk{\sup_{0\le s\le t<\IF}Y(s,t)>u^{1-H}}},  
\EQN
and
\BQN\label{eq:tau11}
\pk{\frac{\tau_2(u)-\tilde{t}_0u}{A(u)}\ge x\Big| \tauu<\IF}
&=&\frac{\pk{\sup_{ t\ge\tilde{t}_0 +xA(u)u^{-1}, s\in[0,t]}Y(s,t)>u^{1-H}}}{\pk{\sup_{0\le s\le t<\IF}Y(s,t)>u^{1-H}}}.  
\EQN
}
In order to complete the proof of \netheo{Thmfl} we need to know the tail asymptotic behaviour of the supremum of the Gaussian random field $Y$ over a region which might depend on $u$. Therefore, \Hr{we shall investigate first} the tail asymptotic behaviour of the supremum of certain \rL{non-stationary} Gaussian random fields (including $Y$ as a special case) over a region depending on $u$ in \netheo{Thmx} followed then by the proof of  \netheo{Thmfl}.

Hereafter, we assume that all \rP{considered Gaussian random fields (or processes)} have almost surely continuous sample paths.
We need to introduce some more notation starting with the well-known Pickands constant $\mathcal{H}_{\alpha}$ \aH{given} by
$$\mathcal{H}_{\alpha}:=\lim_{T\rightarrow\infty} \frac{1}{T}\mathcal{H}_{\alpha}[0,T],\ \ \alpha\in (0,2],$$
where $$\mathcal{H}_{\alpha}[0,T]=\E{ \exp\biggl(\sup_{t\in[0,T]}\Bigl(\sqrt{2}B_{\alpha}(t)-t^{\alpha}\Bigr)\biggr)}\in(0,\IF),\quad
\ T\in(0,\IF),$$
with  $\{B_{\alpha}(t),t\ge0\}$  a  fBm  with Hurst index $\alpha/2\in(0,1]$.
 It is known that $\mlH_1=1$ and $\mlH_2={1}/{\sqrt{\pi}}$, see  Pickands (1969), Albin  (1990), Piterbarg (1996), D\c{e}bicki (2002),
 Debicki et al.\ (2004), \cJI{Mandjes (2007)}, D\c{e}bicki and Mandjes (2011), \xx{Dieker and Yakir (2013)} for various properties of Pickands constant and its generalizations.
Next we introduce another constant, usually referred to as Piterbarg constant, given by
\BQNY
\mathcal{P}_\alpha^a:=\underset{S\rw\IF}\lim\mathcal{P}_\alpha^a[0,S],\quad \alpha\in (0,2],\ a>0, 
\EQNY
where
\BQNY
\mathcal{P}_\alpha^a[S,T]=\E{ \exp\biggl(\sup_{t\in[S,T]}\Bigl(\sqrt{2}B_{\alpha}(t)-(1+a)\abs{t}^{\alpha}\Bigr)\biggr)}\in (0,\IF),\quad
 S<T.
 \EQNY
 It is also known that
\BQN
 \mathcal{P}_1^a=1+\frac{1}{a}\ \ \text{and}\ \ \mathcal{P}_2^a=\frac{1}{2}\left(1+\sqrt{1+\frac{1}{a}}\right)
 \label{eqpp}
\EQN
see e.g., D\c{e}bicki and Mandjes (2003) and D\c{e}bicki and Tabi\'{s} (2011).  As it will be seen in \netheo{Thmx} below both Pickands and Piterbarg constants are important for our study. We  denote by $\Phi(\cdot)$ the standard normal \xH{distribution (of a $N(0,1)$ random variable)}, and further
set  $\Psi(\cdot):=1-\Phi(\cdot)$.
\COM{
In fact, for  positive constants  $a_i,b_i,\alpha_i,\beta_i, i=1,2$, we shall use the following notation, \ccL{with $\Gamma(\cdot)$ the Euler Gamma function.}
\BQNY
\FFab^1=\left\{
 \begin{array}{cc}
 a_1^{1/\alpha_1} \mlH_{\alpha_1} b_1^{-\frac{1}{\beta_1}} \Gamma\LT(\frac{1}{\beta_1}+1\RT),    &\alpha_1<\beta_1,\\
 \piter_{\alpha_1}^{b_1/a_1},    &\alpha_1=\beta_1,\\
 1,    &\alpha_1>\beta_1,
 \end{array}
  \right.\ \
  \FFab^2=\left\{
 \begin{array}{cc}
 a_2^{1/\alpha_2} \mlH_{\alpha_2},    &\alpha_2<\beta_2,\\
 1,    &\alpha_2=\beta_2,
 \end{array}
  \right.
\EQNY
and
\BQNY
\Lambda_1 (x)   =\left\{
 \begin{array}{cc}
 \int_{-\IF}^xe^{-b_2\abs{t}^{\beta_2}}dt,    &\alpha_2<\beta_2,\\
 \piter_{\alpha_2}^{b_2/a_2}(-\IF,x],    &\alpha_2=\beta_2, 
 \end{array}
  \right.\ \
  \Lambda_2 (x)  =\left\{
 \begin{array}{cc}
 \int_x^{\IF}e^{-b_2\abs{t}^{\beta_2}}dt,    &\alpha_2<\beta_2,\\
 \piter_{\alpha_2}^{b_2/a_2}[x,\IF),    &\alpha_2=\beta_2, 
 \end{array}
  \right.\ \ \ x\in\R.
\EQNY
}

In the following we investigate the tail asymptotic behaviour of the supremum of  \rL{non-stationary} Gaussian random fields over a region which is depend on $u$. Our next result is of interest on its own, and furthermore is the key to the proof of Theorem 1.1.

\BT\label{Thmx}
Let $S,T$ be two positive constants, and let  $\{X(s,t),(s,t)\in\lbrack0,S]\times\lbrack0,T]\}$ be a centered \cJ{Gaussian random field}, with standard deviation function $\sigma(\cdot,\cdot)$ and correlation function $r(\cdot,\cdot,\cdot,\cdot)$.
Assume that $\sigma(\cdot,\cdot)$ attains its
maximum on  $\lbrack0,S]\times[0,T]$ at the unique point $(0,t_0)$, with $t_0\in(0,T)$, and further
\BQN\label{eq:var}
\sigma(s,t)=1-b_{1} s^{\beta}(1+o(1))-b_{2}\abs{t-t_0}^{2}(1+o(1))- b_3s\abs{t-t_0} (1+o(1))
\EQN
as $(s,t)\rightarrow(0,t_0)$ for some  constants $\beta\in(1,2)$, and $b_i>0, i=1,2$, $b_3\in\R$ satisfying $b_2+b_3/2>0$. 
Suppose further that
\BQN\label{eq:corr}
r(s,s^{\prime},t,t^{\prime})=1-(a_1|s-s^{\prime}|^{\beta}+a_2|t-t^{\prime}|^{\beta}
)(1+o(1)) \text{\ \ \ as }(s,t),(s^{\prime
},t^{\prime})\rightarrow(0,{t_0})
\EQN
for some constants $a_i>0$, $i=1,2.$ 
\COM{
If  there exist some $\rho>0$ small enough, $\mu\in(0,2]$, and $\CC>0$ such that for any $(s,t), (s',t')\in[0,\rho]\times[t_0-\rho,t_0+\rho]$
\BQNY
\E{X(s,t)-X(s',t')}^2&\le& \CC (\abs{t-t'}^{\mu}+\abs{s-s'}^{\mu}),
\EQNY
then,}
Then, for any $x\in\R$ 
\BQN
\pk{\sup_{(s,t)\in\widetilde{\Delta_x^1}(u)}X(s,t)>u}&=&  \sqrt{\frac{\pi}{b_2}} a_2^{\frac{1}{\beta}}\piter_{\beta}^{b_1/a_1} \H_\beta u^{ \frac{2}{\beta}-1 }  \Psi(u) \Phi(\sqrt{2b_2}x)\oo\label{eq:mainx1}\\
\pk{\sup_{(s,t)\in\widetilde{\Delta_x^2}(u)}X(s,t)>u}&=&  \sqrt{\frac{\pi}{b_2}} a_2^{\frac{1}{\beta}}\piter_{\beta}^{b_1/a_1} \H_\beta u^{ \frac{2}{\beta}-1 }  \Psi(u) \Psi(\sqrt{2b_2}x)\oo\label{eq:mainx2}
\EQN
as $u\rw\IF$, where  $\delta_1(u)=\LT( {\ln u}/{u}\RT)^{\frac{2}{\beta}}$, $\delta_2(u)= {\ln u}/{u}$ and
\BQN\label{eq:Delt}
\Dux=[0,\delu]\times[t_0-\deluu,t_0+xu^{-1}], \quad  \Duux=[0,\delu]\times[t_0+xu^{-1},t_0+\deluu].
\EQN
\ET
\begin{remarks}  a) If $\beta\in(0,1)$, then \eqref{eq:var} becomes
\BQN
\sigma(s,t)=1-b_{1} s^{\beta}(1+o(1))-b_{2}\abs{t-t_0}^{2}(1+o(1)) \ \ \rP{\text{as}\ (s,t)\to(0,t_0)}. 
\EQN
We mention that in this case both \eqref{eq:mainx1} and \eqref{eq:mainx2} are still valid.

b) It can be shown  along the proof of \netheo{Thmx} that if $x=x(u)$ satisfies the following two conditions
\BQN\label{xu}
\lim_{u\rw\IF}x(u)=\IF,\ \ \ \ x(u)=o(u^{\epsilon}) \ \asu,\ \text{for any}\  \epsilon>0,
\EQN
then \eqref{eq:mainx1} still holds with $\Phi(\sqrt{2b_2}x)$ replaced by 1. Similarly, if $x=-x(u)$ with $x(u)$ satisfying \eqref{xu},  then \eqref{eq:mainx2}  holds  with $\Psi(\sqrt{2b_2}x)$ replaced by 1.
\end{remarks}

\bigskip

\prooftheo{Thmfl} \ccL{Define }
\BQNY
T_1(u)=\inf\{t\ge0:\ Y(s,t)>u^{1-H}\ \text{for some} \  s\in[0, t]\}\ 
\EQNY
and
\BQNY
 T_2(u)=\sup\{t\ge0: Y(s,t)>u^{1-H}\ \text{for some} \  s\in[0, t]\}.
\EQNY
Clearly $\tau_i(u)\overset{d}=u T_i(u), i=1,2$, \rP{with $\overset{d}=$ denoting equivalence in distribution.} Consider first the approximation of $\tauu$. For any \ccL{$x\inr$} \rH{and $u>0$} we have
\BQNY
\pk{\frac{\tau_1(u)-\tilde{t}_0u}{A(u)}\le x\Big| \tauu<\IF}&=&
\rH{\pk{T_1(u)\le \tilde{t}_0 + xA(u) u^{-1}   \Big| T_1(u)<\IF}}\\
&=&\frac{\pk{\sup_{0\le s\le t\le\tilde{t}_0  +xA(u)u^{-1}}Y(s,t)>u^{1-H}}}{\pk{\ccL{\tau_1(u)}<\IF}}.
\EQNY
In view of Hashorva et al.\ (2013) for any $H, \Ga \in(0,1)$
\BQN\label{th2}
\pk{\tauu<\IF}=\pk{\sup_{t\ge 0}W_\Ga(t)>u}=\mathcal{W}_{H}(u)\Psi\left(\frac{c^H u^{1-H}}{H^H (1-H)^{1-H}}\right)(1+o(1))\ \ \text{as}\  u\rw\infty,
\EQN
where
$$
\mathcal{W}_{H}(u)=
 2^{\frac{1}{2}-\frac{1}{2H}}\frac{\sqrt{\pi}}{\sqrt{H(1-H)}} \mathcal{H}_{{2H}}\piter_{2H}^{\frac{1-\Ga}{\Ga}}\left(\frac{c^H u^{1-H}}{H^H (1-H)^{1-H}}\right)^{(1/H-1)}.
$$
Next, we focus on the analysis of $\pk{\sup_{0\le s\le t\le\tilde{t}_0  +xA(u)u^{-1}}Y(s,t)>u^{1-H}}$. By Bonforroni's inequality
\BQN\label{eq:p123}
p_3(u)\le\pk{\sup_{0\le s\le t\le\tilde{t}_0  +xA(u)u^{-1}}Y(s,t)>u^{1-H}}\le p_1(u)+p_2(u)+p_3(u),
\EQN
where $p_i(u), i=1,2,3,$ are defined in \eqref{Borell1}, \eqref{PP} and \eqref{p3u} below. In the following, we shall give the asymptotics of $p_3(u)$ as $u\rw\IF$, and give bounds for both $p_1(u)$ and $p_2(u)$ for $u$ large, assuring that they are relatively negligible.

We first consider bounds for $p_1(u)$ and $p_2(u)$. Since \rL{on the set $\{(s,t): 0\le s\le t<\IF\}$ the  maximum  of the variance function $V_Y^2(s,t)$ is attained uniquely} at $(0,\tilt)$, we obtain from the
Borell-TIS inequality (e.g., Adler and Taylor (2007)) that for any constant $K\ge2\tilt$,
there exist  constants $\rho>0$ small enough and $\theta\in(0,1)$ such that, for $u$ sufficiently large
\BQN\label{Borell1}
p_1(u):=\pk{\underset{s\in[\rho, K]\ \text{or}\ t\in[0,\tilt-\rho]}{\sup_{0\le s\le t\le K}}Y(s,t)>u^{1-H}}\le \exp\LT(-\frac{(u^{1-H}-d)^2}{2\theta V_Y^2(0,\tilt)}\RT),
\EQN
with $d=\E{\sup_{0\le s\le t\le K} Y(s,t)}<\IF.$
If follows that
 \BQN\label{P1}
1-\frac{V_Y(s,t)}{V_Y(0, \tilde{t}_0)}=\left\{
              \begin{array}{ll}
\frac{c^2(1-H)^3}{2H}(\tildet-t)^2(1+o(1))+\frac{(\Ga-\Ga^2)(1-H)^{2H}c^{2H}}{2H^{2H}} s^{2H}(1+o(1))    , &  H\le1/2 ,\\
\frac{c^2(1-H)^3}{2H}(\tildet-t+\Ga s)^2(1+o(1))+\frac{(\Ga-\Ga^2)(1-H)^{2H}c^{2H}}{2H^{2H}} s^{2H}(1+o(1)), & H>1/2
              \end{array}
            \right.
\EQN
as $(s,t)\rw (0,\tildet)$ and further the correlation function of $Y$ satisfies
\BQN\label{P2}
1-Cov\LT(\frac{Y(s,t)}{V_Y(s,t)},\frac{Y(s',t')}{V_Y(s',t')}\RT)=\frac{1}{2\tildet^{2H}}\left(\mid t-t'\mid^{2H}+\Ga^{2}\mid s-s'\mid^{2H}\right)(1+o(1))
\EQN
as $(s,t), (s',t')\rw (0,\tildet)$. \rL{In addition, for the chosen $\rho>0$ small enough there exists some $\CC>0$ such that for any $(s,t), (s',t')\in[0,\rho]\times[\tilt-\rho,\tilt+\rho]$
\BQN\label{P3}
\E{Y(s,t)-Y(s',t')}^2&\le& \CC (\abs{t-t'}^{2H}+\abs{s-s'}^{2H}).
\EQN
}
Next, let
$$A=\frac{H^{1/2}}{c(1-H)^{3/2}},\ \ \  \tilu=\frac{u^{1-H}}{V_Y(0,\tilt)}.$$
 \rL{In the light of \eqref{P1} and \eqref{P3}, by the Piterbarg inequality (see Theorem 8.1 in Piterbarg (1996) or  Theorem 8.1 in Piterbarg (2001))}  
 for all $u$ sufficiently large
\BQN\label{PP}
p_2(u):=\pk{\underset{s\in[ \tilde{\delta}_1(\tilu),\rho]\ \text{or}\ t\in[\tilt-\rho,\tilt-\tilde{\delta}_2(\tilu)]}{\sup_{(s,t)\in[0,\rho]\times[\tilt-\rho,\tilt+xA(u)u^{-1}]}}Y(s,t)>u^{1-H}} \le C_1 \ u^{\frac{2(1-H)}{H}} \exp\left(-\frac{u^{2(1-H)}}{2V_Y^2(0,\tilt)}-C_2(\ln u)^2\right)
\EQN
for some positive constants $C_i,i=1,2$, \rP{where $\tilde{\delta}_1(\tilu)=(\ln \tilu/\tilu)^{1/H}, \tilde{\delta}_2(\tilu)= \ln \tilu/\tilu$.} 
\rP{Further, we have}
\BQN\label{p3u}
p_3(u):=\pk{\sup_{(s,t)\in[0,\tilde{\delta}_1(\tilu)]\times[\tilt-\tilde{\delta}_2(\tilu),\tilde{t}_0   +xA(u)u^{-1}]}Y(s,t)>u^{1-H}}=\pk{\sup_{(s,t)\in \widehat{\Delta_{Ax}^1}(\tilu)}\frac{Y(s,t)}{V_Y(0,\tilt)}>\tilu},
\EQN
where
\rP{$\widehat{\Delta_{Ax}^1}(\tilu)=[0,\tilde{\delta}_1(\tilu)]\times[\tilt-\tilde{\delta}_2(\tilu),\tilt+Ax \tilu^{-1}].$} 
\rL{Utilizing \eqref{P1} and \eqref{P2}} we obtain from \netheo{Thmx} that
\BQN\label{eq:tilu}
\pk{\sup_{(s,t)\in \widehat{\Delta_{Ax}^1}(\tilu)}\frac{Y(s,t)}{V_Y(0,\tilt)}>\tilu}
&=&\mathcal{H}_{2H}\mathcal{P}_{2H}^{\frac{1-\Ga}{\Ga}} 2^{-\frac{1}{2H}}\sqrt{2\pi}A
\frac{c(1-H)}{H} \Psi(\tilu)   \tilu^{\frac{1}{H}-1}\Phi(x)  \oo 
\EQN
as $u\rw\IF$. Consequently, we conclude from \eqref{eq:p123}-\eqref{Borell1}, \eqref{PP}-\eqref{eq:tilu} that
\BQNY
\pk{\sup_{0\le s\le t\le\tilde{t}_0  +xA(u)u^{-1}}Y(s,t)>u^{1-H}}=\mathcal{H}_{2H} \mathcal{P}_{2H}^{\frac{1-\Ga}{\Ga}}2^{-\frac{1}{2H}}\sqrt{2\pi}A
\frac{c(1-H)}{H} \Psi(\tilu)   \tilu^{\frac{1}{H}-1}\Phi(x)  \oo
\EQNY
as $u\rw\IF$, and thus in the light of \eqref{th2}
\BQNY
\rH{\limit{u} \sup_{x\inr}} \ABs{ \pk{\frac{\tau_1(u)-\tilde{t}_0u}{A(u)}\le x\Big| \tauu<\IF}- \Phi(x)}=0.
\EQNY
Using similar arguments, we conclude by the properties of the random field $Y$ and \eqref{eq:mainx2} that
\BQNY
\pk{\sup_{t\ge\tilde{t}_0  +xA(u)u^{-1},s\in[0,t]}Y(s,t)>u^{1-H}}=\mathcal{H}_{2H} \mathcal{P}_{2H}^{\frac{1-\Ga}{\Ga}} 2^{-\frac{1}{2H}}\sqrt{2\pi}A
\frac{c(1-H)}{H} \Psi(\tilu)   \tilu^{\frac{1}{H}-1}\Psi(x)  \oo
\EQNY
as $u\rw\IF$, where we used the fact that for any large enough integer $K>\tilde{t}_0$
\BQNY
\pk{\sup_{0\le s \le t<\IF} Y(s,t)>u^{1-H}}=\pk{\sup_{0\le s \le t<K} Y(s,t)>u^{1-H}}\oo\ \ \ \asu,
\EQNY
see Hashorva et al. (2013).
 \ccL{Therefore}
\BQNY
\pk{\frac{\tau_2(u)-\tilde{t}_0u}{A(u)}\le x\Big| \tauu<\IF}&=&1- \pk{\frac{\tauuu-\tilde{t}_0u}{A(u)}\ge x\Big| \tauu<\IF}\\
&=& \rH{ 1- \pk{T_2(u) \ge \tilde{t}_0 + x A(u)u^{-1}  \Big| T_1(u) <\IF}}\\
&=&1- \frac{\pk{\sup_{ t\ge\tilde{t}_0 +xA(u)u^{-1}, s\in[0,t]}Y(s,t)>u^{1-H}}}{\pk{\ccL{\tau_1(u)}<\IF}}\\  
&\rw & \Phi(x)\ \ \ \asu
\EQNY
 for any $x\inr$. \rH{Hence the proof follows by a direct application of \nelem{lemJoint} below.} \QED

\rH{\BEL \label{lemJoint}
Let $(Z_{u1},Z_{u2}), u>0$ be a bivariate random sequence such that $Z_{u2} \ge Z_{u1} $ almost surely for all large $u$. If the following convergence in distribution
$$ Z_{ui} \todis \mathcal{Z} \quad \ \asu$$
holds for $i=1,2$ with $\mathcal{Z}$ a non-degenerate random variable, then we have the joint convergence in distribution
\BQN
(Z_{u1}, Z_{u2}) \todis (\mathcal{Z},\mathcal{Z}) \quad \ \asu.
\EQN
\EEL
}
 {\bf Proof}: \ccL{Let $x,y$ be any two continuous points of the distribution function   $\pk{\mathcal{Z} \le t},  t\inr$.
 It is sufficient to show that
 \BQNY
\rH{\limit{u}} \pk{ Z_{u1} \le x, Z_{u2} \le y}= \pk{\mathcal{Z}\le \min(x,y)}.
 \EQNY}
 In fact, if $x\ge y$ 
by the assumption that $Z_{u2}\ge Z_{u1}$ holds for all large $u$ we have 
\BQNY
\pk{ Z_{u1} \le x, Z_{u2} \le y} 
& =& \pk{Z_{u2} \le  y } \to  \pk{\mathcal{Z} \le y}\ \ \ \ \asu.
\EQNY
 Further, if $x\le y$ 
\BQNY
\pk{ Z_{u1} \le x, Z_{u2} \le y} &=& \pk{Z_{u1} \le x}- \pk{ Z_{u1} \le x, Z_{u2} > y}\\
 &\ge& \pk{Z_{u1} \le x}- \pk{ Z_{u1} \le y, Z_{u2} > y}\\
 &=& \pk{Z_{u1} \le x}- \Bigl( \pk{ Z_{u2} > y} - \pk{ Z_{u1} > y}\Bigr)\\
& \to & \pk{\mathcal{Z} \le  x } \ \ \ \ \asu
\EQNY
and
\BQNY
\pk{ Z_{u1} \le x, Z_{u2} \le y}&\le &\pk{Z_{u1} \le x}  \to  \pk{\mathcal{Z} \le  x } \ \ \ \ \asu
\EQNY
hold, hence the claim follows.
\QED

\COM{

The next theorem  gives asymptotics of the supremum of the mean-value Gaussian random fields with a variance function having unique maximum point on a compact set, where the variance function does not possess a $(E,\alpha)$ structure as introduced in Piterbarg (1996), and thus the theorem complements the result of Piterbarg (1996).

\BT\label{ThmGPiter}
Let $\{X(s,t),(s,t)\in\lbrack0,S]\times\lbrack0,T]\}$ be a centered Gaussian random field as in \netheo{Thmx} with variance function taking its unique maximum at point $(s_0,t_0)$ \xx{in} $[0,S]\times[0,T]$. Suppose that  (\ref{eq:var})-(\ref{eq:corr}) are satisfied with the parameters therein.
Assume further that
there exist some positive constants $G, \mu$ with $\mu\in(0,2]$, and $\rho$ small enough such that
\BQN\label{eq:Inc}
\E{(X(s,t)-X(s',t'))^2}\le G(|s-s^{\prime}|^{\xx{\mu}} +|t-t^{\prime}|^{\xx{\mu}})
\EQN
for any $(s,t), (s',t')\in D_\rho$ with $D_\rho:=[\max(s_0-\rho,0),\min(s_0+\rho,T)]\times[\max(t_0-\rho,0),\min(t_0+\rho,T)]$.
If
$
{\beta_3}/{\beta_1}+{\beta_4}/{\beta_2}>1
$ and one of the following two conditions holds when $b_3<0$:

C1)  $\beta_1=2\beta_3$ and $b_1+b_3/2>0,$\\
 C2)  $\beta_2=2\beta_4$ and $b_2+b_3/2>0$,

then
\BQN\label{eq:main}
\pk{\sup_{(s,t)\in[0,S]\times[0,T]}X(s,t)>u}= \prod_{i=1}^2\LT(\FFab^i u^{\LT(\frac{2}{\alpha_i}-\frac{2}{\beta_i}\RT)_+ }\RT)  \Psi(u)\oo,
\EQN
as $u\rw\IF$, where
\BQNY
\FFab^i=\left\{
 \begin{array}{cc}
 a_i^{1/\alpha_i} \mlH_{\alpha_i} \int_{T_i^{\IF}}e^{-b_i\abs{t}^{\beta_i}}dt,    &\alpha_i<\beta_i,\\
 \piter_{\alpha_i}^{b_i/a_i}[T_i^{\IF}],    &\alpha_i=\beta_i,\\
 1,    &\alpha_i>\beta_i,
 \end{array}
  \right.\ \ \ i=1,2,
\EQNY
with $T_1^{\infty}=\underset{{x\rw\infty}}\lim x([0,S]-s_0)$ and $T_2^{\infty}=\underset{{x\rw\infty}}\lim x([0,T]-t_0)$ (e.g., $\underset{x\rw\infty}\lim\ x[0,1]=[0,\IF)$).
\ET

\begin{remarks}\label{b3n} a)
\netheo{ThmGPiter}, which gives sufficient conditions assuring that \eqref{eq:main} does not depend on the value of the constant $b_3$, extends Theorem 2.1 of Hashorva et al. (2013). \\
b) The notation $T_i^\IF, \cL{i=1,2}$ are introduced in order to give a more compact \cL{result} dealing with both boundary and non-
boundary points cases simultaneously.
\end{remarks}

}

\section{Proof of \netheo{Thmx} }


\def\xiS{ \xH{\tilde \xi(s,t)}}
\prooftheo{Thmx}
We present only the proof of \eqref{eq:mainx1} with $x\ge0$, since the \rH{other cases} can be \rH{dealt} with \ccL{using} the same argumentations.
\rH{For simplicity we shall} assume that $a_1=a_2=1$;  the general case \rH{follows} by a time scaling.

\rP{Since our approach is asymptotic in natural and that $\delta_1(u)$ and $\delta_2(u)$   both converge to 0 as $u$ tends to infinity, the properties \eqref{eq:var} and \eqref{eq:corr} are the only necessary properties of the Gaussian random field $X$ needed for the asymptotics (which can be seen from the proof below).
\ccL{Therefore, we conclude that} 
\BQNY
\pk{\sup_{(s,t)\in\Dux}X(s,t)>u}= \pk{\sup_{(s,t)\in\Dux} \xiS>u} (1+o(1))=: \pi(u)  (1+o(1))\ \ \ \ \asu,
\EQNY
with $\{\xiS, s,t\ge0\}$ any Gaussian random field possessing the properties \eqref{eq:var} and \eqref{eq:corr}. Particularly, we set}
$$
 \xiS=\frac{\xi(s,t)}{(1+b_1s^{\beta})(1+b_2\abs{t-t_0}^{2}+b_3\abs{t-t_0}s)}, \ \ s,t\ge0,
 $$
with  $\{\xi(s,t), s,t\ge0\}$ a \rL{centered Gaussian random field} with covariance function
$$
r_\xi(s,t)=\exp(- s^{\beta}- t^{\beta}), \quad s,t\ge0.
$$
Since $\beta<2$, for any positive constants $S_1, S_2$, we can divide the intervals $[0,\delu]$ and $[t_0-\deluu,t_0+xu^{-1}]$ into several sub-intervals of length $S_1u^{-2/\beta}$ and $S_2u^{-2/\beta}$, respectively.  Specifically, let for $S_1,S_2>0$
\def\det{\triangle}
\def\wdet{\widetilde{\triangle}}
\def\H{\mathcal{H}}
\BQNY
\det_0^i=u^{-\frac{2}{\beta}}[0,S_i], \ \ \ \det_k^i=u^{-\frac{2}{\beta}}[k S_i,(k +1)S_i], \ \ k \in\mathbb{Z}, \ i=1,2.
\EQNY
Let further for any $u>0$
\BQNY
h_1(u)=\lfloor S_1^{-1}(\ln u)^{\frac{2}{\beta}}\rfloor+1,\
h_2(u)=\lfloor S_2^{-1} (\ln u)  u^{\frac{2}{\beta}-1}\rfloor+1,\ i=1,2,\  h_{2,x}(u)=\lfloor S_2^{-1} x u^{\frac{2}{\beta}-1}\rfloor+1.
\EQNY
Here $\lfloor\cdot\rfloor$ denotes the ceiling function. Applying Bonferroni's inequality we obtain
\BQNY
\pi(u)&\le& \sum_{k_1=0}^{h_1(u)}\sum_{k_2=-h_2(u)}^{h_{2,x}(u)}\pk{\underset{(s,t)\in\det_{k_1}^1\times
(t_0+\det_{k_2}^2)}\sup\xiS
>u}\\
&=&\sum_{k_2=-h_2(u)}^{h_{2,x}(u)}\pk{\underset{(s,t)\in\det_{0}^1\times(t_0+\det_{k_2}^2)}\sup
\xiS 
>u}
+\sum_{k_1=1}^{h_1(u)}\sum_{k_2=-h_2(u)}^{h_{2,x}(u)}\pk{\underset{(s,t)\in\det_{k_1}^1\times(t_0+\det_{k_2}^2)}\sup
\xiS
>u}\\
&=:&I_{1,x}(u)+I_{2,x}(u)
\EQNY
and
\BQNY
\pi(u)&\ge& \sum_{k_2=-h_2(u)+1}^{h_{2,x}(u)-1}\pk{\underset{(s,t)\in\det_{0}^1\times(t_0+\det_{k_2}^2)}\sup
\xiS
>u}\\
&&- \sum_{-h_2(u)+1\le i<j\le h_{2,x}(u)-1}
\pk{\underset{(s,t)\in\det_{0}^1\times(t_0+\det_{i}^2)}\sup\xiS
>u,\underset{(s,t)\in\det_{0}^1\times(t_0+\det_{j}^2)}\sup\xiS 
>u}\\
&=:&J_{1,x}(u)- J_{2,x}(u).
\EQNY
Next we derive the \cJI{required} asymptotic bounds of
$I_{1,x}(u)$ and $J_{1,x}(u)$, and show that
\BQN\label{eq:double0x}
I_{2,x}(u)=J_{2,x}(u)(1+o(1))=o(I_{1,x}(u))=o(J_{1,x}(u))\ \ \ \text{as}\  u\rw\IF,\ S_i\rw\IF, i=1,2.
\EQN
Assuming further that $b_3>0$, we  have
\BQNY
J_{1,x}(u)
&\ge&\sum_{k_2=0}^{h_{2,x}(u)-1}\pk{\underset{(s,t)\in\det_{0}^1\times\det_{k_2}^2}\sup\frac{\xi(s,t)}{1+b_1s^{\beta}}
>u(1+b_2((k_2+1)S_2u^{-\frac{2}{\beta}})^{2}+b_3((k_2+1)S_2u^{-\frac{2}{\beta}})(S_1u^{-\frac{2}{\beta}}))}\nonumber\\
&&+\sum_{k_2=-h_2(u)+1}^{-1}\pk{\underset{(s,t)\in\det_{0}^1\times\det_{k_2}^2}\sup\frac{\xi(s,t)}{1+b_1s^{\beta}}
>u(1+b_2(-k_2S_2u^{-\frac{2}{\beta}})^{2}+b_3(-k_2S_2u^{-\frac{2}{\beta}})(S_1u^{-\frac{2}{\beta}}))}\nonumber\\
&&=:J_{1,1,x}(u)+J_{1,2,x}(u).
\EQNY
In view of  \nelem{lemma6.1} in Appendix
\BQN\label{eq:J11}
J_{1,1,x}(u)&=& \mathcal{P}_{\beta}^{b_1}[0,S_1]\mathcal{H}_{\beta}[0,S_2]\frac{1}{\sqrt{2\pi} u}\sum_{k_2=0}^{h_{2,x}(u)-1}
\frac{1}{1+b_2((k_2+1)S_2u^{-\frac{2}{\beta}})^{2}+b_3((k_2+1)S_2u^{-\frac{2}{\beta}})(S_1u^{-\frac{2}{\beta}})}\nonumber \\ &&\times\exp\left(-\frac{u^2(1+b_2((k_2+1)S_2u^{-\frac{2}{\beta}})^{2}+b_3((k_2+1)S_2u^{-\frac{2}{\beta}})(S_1u^{-\frac{2}{\beta}}))^2}{2}\right)\oo\nonumber\\
&=& \mathcal{P}_{\beta}^{b_1}[0,S_1]\mathcal{H}_{\beta}[0,S_2] \Psi(u) \nonumber\\     
&&\times \sum_{k_2=0}^{h_{2,x}(u)-1}\exp\left(-b_2((k_2+1)S_2u^{1-\frac{2}{\beta}})^{2}-
b_3u^2((k_2+1)S_2u^{-\frac{2}{\beta}})(S_1u^{-\frac{2}{\beta}})\right)\oo\nonumber\\
&=&\mathcal{P}_{\beta}^{b_1}[0,S_1]\frac{\mathcal{H}_{\beta}[0,S_2]}{S_2} \Psi(u)   u^{\frac{2}{\beta}-1} \int_0^x e^{-b_2 y^{2}}dy\oo
\EQN
as $u\rw\IF,$
where in the last equation we utilised the facts that
\BQNY\label{eqkey}
h_{2,x}(u)\rw\IF,\ \ h_{2,x}(u) S_2u^{1-\frac{2}{\beta}}\rw x,\ \ u^2( h_{2,x}(u) S_2u^{-\frac{2}{\beta}})(S_1u^{-\frac{2}{\beta}}) \rw 0
\EQNY
as $u\rw\IF$. Similarly
\BQN\label{eq:J12}
J_{1,2,x}(u)
&=&\mathcal{P}_{\beta}^{b_1}[0,S_1]\frac{\mathcal{H}_{\beta}[0,S_2]}{S_2} \Psi(u)   u^{\frac{2}{\beta}-1} \int_{-\IF}^0  e^{-b_2y^{2}}dy\oo
\EQN
as $u\rw\IF$. Therefore we conclude that
\BQN\label{eq:J_1x}
J_{1,x}(u)
&\ge&\mathcal{P}_{\beta}^{b_1}[0,S_1]\frac{\mathcal{H}_{\beta}[0,S_2]}{S_2} \Psi(u)   u^{\frac{2}{\beta}-1} \int_{-\IF}^x  e^{-b_2y^{2}}dy\oo\ \ \asu.
\EQN
Using similar arguments we further obtain that
\BQN\label{eq:I_1x}
I_{1,x}(u)
&\le&\sum_{k_2=0}^{h_{2,x}(u)-1}\pk{\underset{(s,t)\in\det_{0}^1\times\det_{k_2}^2}\sup\frac{\xi(s,t)}{1+b_1s^{\beta}}
>u(1+b_2(k_2S_2u^{-\frac{2}{\beta}})^{2})}\nonumber\\
&&+\sum_{k_2=-h_2(u)}^{-1}\pk{\underset{(s,t)\in\det_{0}^1\times\det_{k_2}^2}\sup\frac{\xi(s,t)}{1+b_1s^{\beta}}
>u(1+b_2(-(k_2+1)S_2u^{-\frac{2}{\beta}})^{2})}\nonumber\\
&=&\mathcal{P}_{\beta}^{b_1}[0,S_1]\frac{\mathcal{H}_{\beta}[0,S_2]}{S_2} \Psi(u)   u^{\frac{2}{\beta}-1} \int_{-\IF}^x  e^{-b_2y^{2}}dy\oo
\EQN as $u\rw\IF$.
Next we verify \eqref{eq:double0x}. Specifically
\BQNY
I_{2,x}(u)
&\le&\sum_{k_1=1}^{h_1(u)}\sum_{k_2=0}^{h_{2,x}(u)}\pk{\underset{(s,t)\in\det_{k_1}^1\times\det_{k_2}^2}\sup \xi(s,t) >u(1+b_1(k_1S_1u^{-\frac{2}{\beta}})^{\beta}+b_2(k_2S_2u^{-\frac{2}{\beta}})^{2})}\\
&&+\sum_{k_1=1}^{h_1(u)}\sum_{k_2=-h_2(u)}^{-1}\pk{\underset{(s,t)\in\det_{k_1}^1\times\det_{k_2}^2}\sup \xi(s,t) >u(1+b_1(k_1S_1u^{-\frac{2}{\beta}})^{\beta}+b_2(-(k_2+1)S_2u^{-\frac{2}{\beta}})^{2})}. 
\EQNY
Similar argumentations as in \eqref{eq:J_1x} yield
\BQN\label{I2x}
I_{2,x}(u)\le\H_{\beta}[0,S_1]\H_{\beta}[0,S_2]\Psi(u)(S_2^{-1}u^{\frac{2}{\beta}-1})\int_{-\IF}^x  e^{-b_2y^{2}}dy\sum_{k_1=1}^{h_1(u)}
\exp\left(-b_1 (k_1S_1)^{\beta}\right)\oo
\EQN
as $u\rw\IF$.
Further, we write
\BQNY
J_{2,x}(u)&=&\sum_{-h_2(u)+1\le i<j\le h_{2,x}(u)-1}\pk{\underset{(s,t)\in\det_{0}^1\times(t_0+\det_{i}^2)}\sup\xiS
>u,
\underset{(s,t)\in\det_{0}^1\times(t_0+\det_{j}^2)}\sup\xiS
>u}\\
&=:&\Sigma_{1,x}(u)+\Sigma_{2,x}(u),
\EQNY
where $\Sigma_{1,x}(u)$ is the sum over indexes $j=i+1$, and $\Sigma_{2,x}(u)$ is the sum over indexes $j>i+1$.
Let
$$
B(i,S_2,u)=u(1+b_2(\abs{i}S_2u^{-\frac{2}{\beta}})^{2}), \ \ i\in \mathbb{Z},\ S_2>0,\ u>0.
$$
It follows that
\BQNY
\Sigma_{1,x}(u)&\le&\sum_{i=-1}^{ h_{2,x}(u)-1}\pk{\underset{(s,t)\in\det_{0}^1\times\det_{i}^2}\sup\frac{\xi(s,t)}{1+b_1s^{\beta}}>B(0,S_2,u),
\underset{(s,t)\in\det_{0}^1\times\det_{i+1}^2}\sup\frac{\xi(s,t)}{1+b_1s^{\beta}}>B(0,S_2,u)}\\
&&+\sum_{i=- h_2(u)+1}^{ -2}\pk{\underset{(s,t)\in\det_{0}^1\times\det_{i}^2}\sup\frac{\xi(s,t)}{1+b_1s^{\beta}}>B(i+2,S_2,u),
\underset{(s,t)\in\det_{0}^1\times\det_{i+1}^2}\sup\frac{\xi(s,t)}{1+b_1s^{\beta}}>B(i+2,S_2,u)}
\EQNY
and, for any $i,j\in \mathbb{Z}$
\BQNY
&&\pk{\underset{(s,t)\in\det_{0}^1\times\det_{i}^2}\sup\frac{\xi(s,t)}{1+b_1s^{\beta}}>B(j,S_2,u),
\underset{(s,t)\in\det_{0}^1\times\det_{i+1}^2}\sup\frac{\xi(s,t)}{1+b_1s^{\beta}}>B(j,S_2,u)}\\
&&=\pk{\underset{(s,t)\in\det_{0}^1\times\det_{0}^2}\sup\frac{\xi(s,t)}{1+b_1s^{\beta}}>B(j,S_2,u)}+\pk{
\underset{(s,t)\in\det_{0}^1\times\det_{1}^2}\sup\frac{\xi(s,t)}{1+b_1s^{\beta}}>B(j,S_2,u)}\\
&&-\pk{\underset{(s,t)\in\det_{0}^1\times(\det_{0}^2\cup\det_1^2)}\sup\frac{\xi(s,t)}{1+b_1s^{\beta}}>B(j,S_2,u)}.
\EQNY
Therefore, analogous to the derivation of \eqref{eq:J_1x}, we obtain
\BQN\label{eq:sig1x}
\limsup_{u\rw\IF}\frac{\Sigma_{1,x}(u)}{\Psi(u)   u^{\frac{2}{\beta}-1}}&\le&
\mathcal{P}_{\beta}^{b_1}[0,S_1] \frac{2\mathcal{H}_{\beta}[0,S_2]-\H_{\beta}[0,2S_2]}{S_2}\LT( x+ \int_{-\IF}^0 e^{-b_2y^{2}}dy\RT).
\EQN
Further, for any $u>0$
\BQNY
\Sigma_{2,x}(u)&\le&\sum_{i=-1}^{ h_{2,x}(u)-1}\sum_{j\ge2}\pk{\underset{(s,t)\in\det_{0}^1\times\det_{0}^2}\sup \xi(s,t) >u,
\underset{(s,t)\in\det_{0}^1\times\det_{j}^2}\sup \xi(s,t) >u}\\
&&+\sum_{i=-h_2(u)+1}^{-2}\sum_{j\ge2}\pk{\underset{(s,t)\in\det_{0}^1\times\det_{0}^2}\sup \xi(s,t) >B(i+1,S_2,u),
\underset{(s,t)\in\det_{0}^1\times\det_{j}^2}\sup \xi(s,t) >u}\\
&\le&\sum_{i=- 1}^{ h_{2,x}(u)-1}\sum_{j\ge2}\pk{\underset{(s',t')\in\det_{0}^1\times\det_{j}^2}{\sup_{(s,t)\in\det_{0}^1\times\det_{0}^2}} \zeta(s,t,s',t') >2u}\\
&&+\sum_{i=-h_2(u)+1}^{ -2}\sum_{j\ge2}\pk{\underset{(s',t')\in\det_{0}^1\times\det_{j}^2}{\sup_{(s,t)\in\det_{0}^1\times\det_{0}^2}} \zeta(s,t,s',t') >B(i+1,S_2,u)+u},
\EQNY
where
$$
\zeta(s,t,s',t')=\xi(s,t)+\xi(s',t'),\ \ s,s',t,t'\ge0.
$$
It is easy to check that, for $u$ sufficiently large
$$
2\le\E{(\zeta(s,t,s',t'))^2}=4-2(1-r(\abs{s-s'},\abs{t-t'}))\le 4- ((j-1)S_2)^{\beta} u^{-2}
$$
for any $ (s,t)\in\det_{0}^1\times\det_{0}^2, (s',t')\in\det_{0}^1\times\det_{j}^2$.
Borrowing the arguments of the proof of Lemma 6.3 in Piterbarg (1996) we conclude that
\BQN\label{eq:sig2x}
\limsup_{u\rw\IF}\frac{\Sigma_{2,x}(u)}{\Psi(u)   u^{\frac{2}{\beta}-1}}&\le& \CC\ x\
(\mathcal{H}_{\beta}[0,S_1])^2  S_2 \sum_{j\ge1}\exp\left(-\frac{1}{8}(jS_2)^{\beta}\right)
\EQN for some positive constant $\CC$. Hence  the claim 
follows from  (\ref{eq:double0x}--\ref{eq:sig2x}) when  $b_3>0$ by letting $S_2,S_1\rw\IF$. When
$b_3<0$, the same results can be obtained using similar arguments as above \ccL{and} the  fact that
\BQNY
1-\sigma(s,t)\ge
b_{1}s^{\beta}(1+o(1))+\left(b_{2}+\frac{b_3}{2}\right)\abs{t-t_0}^{2}(1+o(1))
\EQNY
as $(s,t)\rightarrow(0,t_0)$ which is utilised for verifying \eqref{eq:double0x}, and thus the proof is complete.
 \QED

\section{Appendix: Piterbarg Lemma for Gaussian Random Fields}
In order to find the asymptotics of supremum of centered non-smooth Gaussian processes two crucial results are important,
 namely \rL{the Pickands Lemma and the Piterbarg Lemma.}  Although for experts in this field the results are well-known, \rL{we would like to briefly mention them.
Let $\{X(t),t\ge 0\}$ be a} centered stationary Gaussian process
 with a.s.\ continuous sample paths and correlation function $r(t)$ which satisfies $r(t) = 1 -  \cJ{t^{\alpha}}(1+ o(1))$ as $t \rightarrow 0$
 with $\alpha \in (0, 2]$ and $r(t)<1$ for all $t>0$.  In the seminal paper Pickands (1969) it was shown that
 for any $T\in (0,\IF)$
\BQN\label{pic}
\pk{\sup_{t\in[0,T]}X(t)>u}=\Ha T u^{\frac{2}{\alpha}}\Psi(u)(1+o(1))\ \ \ \text{as}\ \ u \rw \IF.
\EQN
{The proof of \eqref{pic} strongly relies on Pickands Lemma which says that
\BQN\label{eq:pick1}
\pk{\sup_{t\in[0,u^{-{\frac{2}{\alpha}}}T]}X(t)>u}=\Ha[0, T] \Psi(u)
(1+ o(1)) \ \ \ \text{as}\ \ u \rw \IF.
\EQN
In the seminal contribution Piterbarg (1972) V.I. Piterbarg rigorously proved  \eqref{pic} and then
extended \eqref{eq:pick1} to a result which we refer to as {Piterbarg} Lemma, namely for any constant $b>0$
\BQNY
\pk{\sup_{t\in[0,u^{-{\frac{2}{\alpha}}}T]}\frac{X(t)}{1+bt^\alpha}>u}=\cJ{\mathcal{P}_{\alpha}^b}[0,T]\Psi(u)
(1+ o(1)) \ \ \ \text{as}\ \ u \rw \IF.
\EQNY
\COM{where
\BQNY\label{pick}
\mathcal{P}_{\alpha}^b[0, T]=\E{ \exp\biggl(\sup_{t\in[0,T]}\Bigl(\sqrt{2}B_\alpha(t)-\cJ{(1+b)}t^{\alpha}\Bigr)\biggr)} \in (0,\IF).
\EQNY
The  positive constant (referred to as {\it Piterbarg constant}) given by
\BQNY
\mathcal{P}_\alpha^b=\lim_{T\rw\IF}\mathcal{P}_{\alpha}^b[0, T] \in (0,\IF)
\EQNY
appears naturally when dealing with the extremes of non-stationary Gaussian processes or Gaussian random fields, see e.g., \cite{Pit96} and our main result below. It is known that $\mlH_1=1$, $\mlH_2={1}/{\sqrt{\pi}}$, and
 \BQN
\mathcal{P}_1^b=1+\frac{1}{b},\ \ \ \  \mathcal{P}_2^b=\frac{1}{2}\left(1+\sqrt{1+\frac{1}{b}}\right), \ \ b>0,
 \label{eqpp}
\EQN
}
Our next result is a variant of Piterbarg Lemma \rL{for two-dimensional case}. \Hr{We omit its proof since it follows with exactly the same arguments as that of Lemma 6.1 in Piterbarg (1996).}
\BEL \label{lemma6.1}
Let $\{\xi(s,t), s,t\ge0\}$ be \Hr{a} centered Gaussian random field
\Hr{with covariance function
$$
r_\xi(s,t)=\exp(- s^{\alpha_1}- t^{\alpha_2}),\  \quad s,t\ge0,\ \ \text{with} \ \alpha_1,\alpha_2\in(0,2].
$$
}
\Hr{Let further $S,T_1,T_2$ be three constants such that $S>0$ and $T_1<T_2$}. Then, for any constants $b_1\ge 0, b_2>0,$ and any positive function $g(u), u\ge0$ satisfying $\lim_{u\rw\IF}g(u)/u=1$, we have
\BQN\label{eq:piter6.1}
\pk{
{\sup_{(s,t)\in [0, u^{-\frac{2}{\alpha_1}}S]\times[u^{-\frac{2}{\alpha_2}}T_1,u^{-\frac{2}{\alpha_2}}T_2]}}\frac{\xi(s,t)}{(1+b_1s^{\alpha_1})(1+b_2t^{\alpha_2})}>g(u)}=\mathcal{P}_{\alpha_1}^{b_1}[0,S]
\mathcal{P}_{\alpha_2}^{b_2}[T_1,T_2]\Psi(g(u))(1+o(1)) 
\EQN
 as $u\rw\IF$.
\EEL
\begin{remark}
In the last formula  we identify $\mathcal{P}_{\alpha_1}^{b_1}[0,S]$ to be  $\mathcal{H}_{\alpha_1}[0,S]$  when $b_1=0$.
\end{remark}

\COM{\BEL\label{lemVarY}
The variance function $V^2_Y(s,t)$ given in \eqref{eq:varY} attaints its unique global maximum on the set $\{(s,t): 0\le s\le t<\IF\}$ at $(0, \tilde{t}_0)$, with $\tilde{t}_0=\frac{H}{c(1-H)}$. Moreover
$$
V_Y(0, \tilde{t}_0)=\frac{H^H (1-H)^{1-H}}{c^H }.
$$
\EEL
}

\COM{
\BEL\label{eq:PPP} The Gaussian random field $\{Y(s,t), s,t\ge0\}$ defined in \eqref{Yst} has the following properties.

P1: The standard deviation function of $Y$ satisfies
 \BQN\label{P1}
1-\frac{V_Y(s,t)}{V_Y(0, \tilde{t}_0)}=\left\{
              \begin{array}{ll}
\frac{c^2(1-H)^3}{2H}(\tildet-t)^2(1+o(1))+\frac{(\Ga-\Ga^2)(1-H)^{2H}c^{2H}}{2H^{2H}} s^{2H}(1+o(1))    , &  H\le1/2 ,\\
\frac{c^2(1-H)^3}{2H}(\tildet-t+\Ga s)^2(1+o(1))+\frac{(\Ga-\Ga^2)(1-H)^{2H}c^{2H}}{2H^{2H}} s^{2H}(1+o(1)), & H>1/2
              \end{array}
            \right.
\EQN
as $(s,t)\rw (0,\tildet)$.

P2: The correlation function of $Y$ satisfies
\BQN\label{P2}
1-Cov\LT(\frac{Y(s,t)}{V_Y(s,t)},\frac{Y(s',t')}{V_Y(s',t')}\RT)=\frac{1}{2\tildet^{2H}}\left(\mid t-t'\mid^{2H}+\Ga^{2}\mid s-s'\mid^{2H}\right)(1+o(1))
\EQN
as $(s,t), (s',t')\rw (0,\tildet)$.

\BQN\label{P3}
\E{Y(s,t)-Y(s',t')}^2&\le& \CC (\abs{t-t'}^{2H}+\abs{s-s'}^{2H}).
\EQN

\EEL
}
\COM{
\BEL\label{LemtailK} Let $\{Y(s,t), s,t\ge0\}$ be defined as in \eqref{Yst}. We have
for any large enough integer $K>\tilde{t}_0$
\BQNY
\pk{\sup_{0\le s \le t<\IF} Y(s,t)>u^{1-H}}=\pk{\sup_{0\le s \le t<K} Y(s,t)>u^{1-H}}\oo\ \ \ \asu.
\EQNY
\EEL

\BT\label{th2} Let  $\{W_\gamma(t),t\ge 0\}$  be defined as in \eqref{Wgam}. We have, for any $H, \Ga \in(0,1)$
\BQN\label{th2}
\pk{\tauu<\IF}=\pk{\sup_{t\ge 0}W_\Ga(t)>u}=\mathcal{W}_{H}(u)\Psi\left(\frac{c^H u^{1-H}}{H^H (1-H)^{1-H}}\right)(1+o(1))\ \ \text{as}\  u\rw\infty,
\EQN
where
$$
\mathcal{W}_{H}(u)=
 2^{\frac{1}{2}-\frac{1}{2H}}\frac{\sqrt{\pi}}{\sqrt{H(1-H)}} \mathcal{H}_{{2H}}\piter_{2H}^{\frac{1-\Ga}{\Ga}}\left(\frac{c^H u^{1-H}}{H^H (1-H)^{1-H}}\right)^{(1/H-1)}.
$$
\ET

\bigskip

}

{\bf Acknowledgement}: The authors kindly acknowledge partial
support from Swiss National Science Foundation Project 200021-1401633/1.
 and the project RARE -318984, a Marie Curie IRSES Fellowship within the 7th European Community Framework Programme.

\end{document}